\documentclass{amsart}

\usepackage{amssymb}
\usepackage{amsmath}
\usepackage{amscd}
\usepackage{multicol}
\usepackage[dvipdfmx]{graphicx}
\usepackage[all]{xy}
\usepackage{tabularx}
\title[Automorphism group computed by Galois points, II]{Automorphism group of plane curve computed by Galois points, II}

\author{Takeshi Harui}
\address{Department of Core Studies, Kochi University of Technology, Kami, Tosayamada, Kochi 782-8502, Japan}
\email{harui.takeshi@kochi-tech.ac.jp}

\author{Kei Miura}
\address{Department of Mathematics, National Institute of Technology, Ube College, Ube, Yamaguchi, 755-8555, Japan}
\email{kmiura@ube-k.ac.jp}
\thanks{The second author was partially supported by JSPS KAKENHI Grant Number JP26400057.}

\author{Akira Ohbuchi}
\address{Department of Mathematical Sciences, Faculty of Science and Technology, Tokushima University, Tokushima, 770-8502, Japan}
\email{ohbuchi@tokushima-u.ac.jp}
\thanks{The third author was partially supported by JSPS KAKENHI Grant Number JP15K04822.}

\date{September 15, 2017}

\keywords{icosahedral group, Galois point, plane curve, automorphism group}
\subjclass[2010]{Primary 14H37; Secondary 14H50}

\newtheorem{theorem}{Theorem}

\newtheorem{proposition}[theorem]{Proposition}

\newtheorem{fact}{Fact}

\theoremstyle{definition}
\newtheorem{definition}{Definition}

\theoremstyle{remark}
\newtheorem{remark}{Remark}

\newenvironment{namelist}[1]{%
\begin{list}{}
  {
   \settowidth{\labelwidth}{#1}
   \setlength{\leftmargin}{2.5\labelwidth}}
}{%
\end{list}}
\def\dfrac#1#2{{\displaystyle\frac{#1}{#2}}}
\begin{document}

\begin{abstract}
Recently, the first author \cite{harui} classified finite groups obtained as automorphism groups of 
smooth plane curves of degree $d \ge 4$ into five types.  
He gave an upper bound of the order of the automorphism group for each types. 
For one of them, the type (a-ii), that is given by {\rm max}$\left\{ 2 d (d - 2), 60 d \right\}$. 
In this article, we shall construct typical examples of smooth plane curve $C$ 
by applying the method of Galois points, whose automorphism group has order $60 d$. 
In fact, we determine the structure of the automorphism group of those curves.  
\end{abstract}

\maketitle

%%%%%%%%%%%%%%%%%%%%%%%%%%%%%%%%%%%%%%%%%%%%%%%%%%%%%%%%%
\section{Introduction}

The purpose of this article is to give typical examples of smooth plane curve of degree $d$ 
whose automorphism group has order $60 d$. 
In fact, we study the structure of that group.  
Our method is based on the classification theorem of automorphism groups by the first author
and the theory of Galois points for smooth plane curves.

First, we recall several definitions of Galois points in brief. 
Throughout the present article, we work over the complex number field ${\mathbb C}$. 
The concept of Galois points was introduced by Yoshihara in 1996 (e.g. \cite{m-y1}). 
Let $C \subset {\mathbb P}^2$ be a smooth plane curve of degree $d$ $(d \ge 4)$ 
and ${\mathbb C} ( C )$ the function field of $C$. 
Let $P$ be a point of ${\mathbb P}^2$. 
Consider the morphism $\pi_P : C \rightarrow {\mathbb P}^1$, 
which is the restriction of the projection 
${\mathbb P}^2 \dashrightarrow {\mathbb P}^1$ with the center $P$. 
Then we obtain the field extension induced by $\pi_P$, 
i.e., $\pi_P^* : {\mathbb C} ({\mathbb P}^1) \hookrightarrow {\mathbb C} (C)$. 
Putting $K_P = \pi_P^* ({\mathbb C} ({\mathbb P}^1))$, we have the following definition.

\begin{definition}
The point $P$ is called a {\em Galois point} for $C$ if the field extension ${\mathbb C} (C) / K_P$ 
is Galois. Furthermore, a Galois point is said to be {\em inner} (resp. {\em outer}) if 
$P \in C$ (resp. $P \in {\mathbb P}^2 \setminus C$). 
The group $G_P = {\rm Gal} ({\mathbb C}(C) / K_P)$ is called the {\em Galois group at $P$}.  
\end{definition}

We denote by $\delta (C)$ (resp. $\delta' (C)$) the number of 
inner (resp. outer) Galois points for $C$. 
There are many known results on Galois points. 
We recall some of them.

\begin{theorem}[\cite{m-y1}, \cite{yoshi}]
Suppose that $C$ is a smooth plane curve of degree $d$ $(d \ge 4)$. Then, 
 \begin{namelist}{11}
  \item[{\rm (i)}] $\delta'(C) = 0, 1$ or $3$. Further, $\delta'(C) = 3$ if and only if $C$ is projectively equivalent to the Fermat curve. 
  \item[{\rm (ii)}] $\delta (C) = 0, 1$ or $4$ if $d = 4$. Further, $\delta(C) = 4$ if and only if $C$ is projectively equivalent to the curve defined by
  $X^4 + Y^4 + Y Z^3 = 0$. When $d \ge 5$, we have $\delta(C) = 0$ or $1$. 
 \end{namelist}
\end{theorem}

\begin{theorem}[\cite{yoshi}]
Suppose that $C$ is a smooth plane curve of degree $d$ $(d \ge 4)$. 
If $P$ is an inner $($resp. outer$)$ Galois point, 
then $G_P$ is isomorphic to the cyclic group of degree $d - 1$ $($resp. $d$$)$, i.e., $G_P \cong {\mathbb Z}_{d - 1}$ 
$($resp. ${\mathbb Z}_d$$)$. 
\end{theorem}

\begin{remark}
If $C$ has singularities, then the theorem above does not hold true. 
Namely, there exist a singular plane curve $C$ and a Galois point $P$ for $C$ 
such that $G_P$ is not cyclic. 
For example, see \cite{miura}. 
\end{remark}

When $C$ has a Galois point, we can give a concrete defining equation of $C$.

\begin{proposition}[\cite{yoshi}]\label{stform}
By a suitable change of coordinates, the defining equation of $C$ with an outer Galois point 
can be expressed as $Z^d + F_d(X, Y) = 0$, where $F_d(X, Y)$ is a homogeneous polynomial 
of degree $d$ without multiple factors. 
\end{proposition}

Referring to \cite{harui}, we may infer that 
plane curves with $\delta (C) \ne 0$ or $\delta' (C) \ne 0$ play an important role 
when we classify the automorphism group of smooth plane curves.

In \cite{harui}, the first author classified finite groups obtained as automorphism groups of $C$ into five types. 
First of all, we recall several definitions. 
Let $G$ be a group of automorphisms of $C$. 
Then, it is well-known that $G$ is considered as a subgroup of ${\rm PGL} (3, {\mathbb C}) = {\rm Aut} ({\mathbb P}^2)$. 
Let $F_d$ be the Fermat curve $X^d + Y^d + Z^d = 0$. 
We denote by $K_d$ a smooth curve defined by $X Y^{d - 1} + Y Z^{d - 1} + Z X^{d - 1} = 0$ 
(In \cite{harui}, $K_d$ is called {\em Klein curve} of degree $d$). 
For a non-zero monomial $c X^i Y^j Z^k$ with $c \in {\mathbb C} \setminus \{ 0 \}$, we define its {\em exponent} as ${\rm max} \left\{ i, j, k \right\}$. 
For a homogeneous polynomial $F(X, Y, Z)$, the {\em core} of $F(X, Y, Z)$ is defined as the sum of all terms of 
$F(X, Y, Z)$ with the greatest exponent.

\begin{definition}
Let $C_0$ be a smooth plane curve with defining equation $F_0 (X, Y, Z) = 0$. 
Then a pair $(C, G)$ of a smooth plane curve $C$ and a subgroup $G \subset {\rm Aut}(C)$ is said to be a 
{\em descendent} of $C_0$ if $C$ is defined by a homogeneous polynomial whose core coincides with $F_0 (X, Y, Z)$ 
and $G$ acts on $C_0$ in a suitable coordinate system. 
\end{definition}

\begin{definition}
We denote by ${\rm PBD}(2, 1)$ the following subgroup of ${\rm PGL} (3, {\mathbb C})$: 
$$
{\rm PBD}(2, 1) := \left. \left\{ 
A = \left(\begin{array}{rrr}
a_{11} &a_{12} &0  \\
a_{21} &a_{22} &0  \\
0      &0      &\alpha
\end{array}\right) \in {\rm GL} (3, {\mathbb C})
\right\} \right/ {\mathbb C}^{\times}.$$
\end{definition}

We remark that there exists a natural group homomorphism 
$\rho  : {\rm PBD}(2, 1) \rightarrow {\rm PGL} (2, {\mathbb C})$, i.e., 
$A \mapsto  (a_{ij})$.

Using these concepts, the first author proved the following theorem.

\begin{theorem}[\cite{harui}]\label{harui_th}
Let $C$ be a smooth plane curve of degree $d \ge 4$, $G$ a subgroup of ${\rm Aut} (C)$. 
Then one of the following holds:
\begin{namelist}{1}
\item[{\rm (a-i)}] $G$ fixes a point on $C$ and $G$ is a cyclic group whose order is at most $d (d - 1)$. 
Furthermore, if $d \ge 5$ and $\vert G \vert = d (d - 1)$, then $C$ is projectively equivalent to the 
curve $Y Z^{d - 1} + X^d + Y^d = 0$ and ${\rm Aut}(C) \cong {\mathbb Z}_{d (d - 1)}$. 
\item[{\rm (a-ii)}] $G$ fixes a point not lying on $C$ and there exists a commutative diagram
$$\xymatrix{
1 \ar[r] &  {\mathbb C}^{\times} \ar[r]  & {\rm PBD}(2, 1) \ar[r]^{\rho} & {\rm PGL} (2, {\mathbb C}) \ar[r] & 1 \\
1 \ar[r] &  N       \ar@{^{(}->}[u] \ar[r] & G        \ar@{^{(}->}[u] \ar[r] & G' \ar@{^{(}->}[u] \ar[r]             & 1, \\
}$$
where $N$ is a cyclic group whose order is a factor of $d$ and $G'$ is a subgroup of ${\rm PGL} (2, {\mathbb C})$,
i.e., a cyclic group ${\mathbb Z}_m$, 
a dihedral group $D_{2 m}$, 
the tetrahedral group $A_4$, 
the octahedral group $S_4$ 
or the icosahedral group $A_5$. 
Furthermore, $m \le d - 1$ and if $G' \cong D_{2 m}$ then $m \mid d - 2$ or $N$ is trivial. 
In particular, $\vert G \vert \le$ ${\rm max} \left\{ 2 d (d - 2), 60 d \right\}$.  
\item[{\rm (b-i)}] $(C, G)$ is a descendant of the Fermat curve $F_d : X^d +Y^d + Z^d = 0$. 
In this case $\vert G \vert \le 6 d^2$. 
\item[{\rm (b-ii)}] $(C, G)$ is a descendant of the Klein curve $K_d : X Y^{d -1} +Y Z^{d - 1} + Z X^{d - 1} = 0$. 
In this case $\vert G \vert \le 3 (d^2 - 3 d + 3)$ if $d \ge 5$. 
On the other hand, $\vert G \vert \le 168$ if $d = 4$.  
\item[{\rm (c)}] $G$ is conjugate to a finite primitive subgroup of ${\rm PGL}(3, {\mathbb C})$. 
Namely, the icosahedral group $A_5$, the Klein group of order $168$, the alternating group $A_6$, 
the Hessian group $H_{216}$ or its subgroup of order $36$ or $72$. 
In particular, $\vert G \vert \le 360$. 
\end{namelist}
\end{theorem}

%%%%%%%%%%%%%%%%%%%%%%%%%%%%%%%%%%%%%%%%%%%%%%%%%%%%%%%%%
\section{Remark on {\rm (a-i)}}

Let $P_1, \cdots , P_m$ be all inner and outer Galois points for $C$ and 
$G(C)$ denote the group generated by $G_{P_i}$ $(i = 1, 2, \ldots, m)$. 
The group $G(C)$ is called the group generated by automorphisms belonging to all Galois points for $C$. 
In \cite{mioh}, we have studied the difference between ${\rm Aut} (C)$ and $G(C)$. 
Referring to \cite{fuka}, 
if $\delta(C) \ge 1$ and $\delta' (C) \ge 1$, then $C$ is projectively equivalent to the curve as in 
Theorem \ref{harui_th} {\rm (a-i)}. We denote the curve by $C(d)$, i.e., $C(d) : Y Z^{d - 1} + X^d + Y^d = 0$. 
If $d \ge 5$, then $P = (0:0:1)$ is the only inner Galois point and 
$Q = (1:0:0)$ is the only outer Galois point for $C(d)$. 
We put $G_P = \left\langle \sigma \right\rangle$ and $G_Q = \left\langle \tau \right\rangle$. 
Then $G(C(d)) = \left\langle \sigma, \tau \right\rangle$. 
In \cite{mioh}, we obtain ${\rm Aut}(C(d)) = G(C(d))$. 
Thus Galois points play an important role in studying of the automorphism groups of smooth plane curves.

%%%%%%%%%%%%%%%%%%%%%%%%%%%%%%%%%%%%%%%%%%%%%%%%%%%%%%%%%
\section{Main results}

In this section, we first remark on Theorem \ref{harui_th} {\rm (a-ii)}. 
In general, we have $2 d (d - 2) > 60 d$. 
However, clearly $2 d (d - 2) < 60 d$ if $d < 32$. 
Hence we consider the case $d < 32$, and try to construct $C$ with $\vert {\rm Aut} (C) \vert = 60 d$.

Let $F_i(X, Y)$ ($i = 1, 2, 3$) be the homogeneous polynomials of $X$ and $Y$ defined by 
\begin{namelist}{11}
\item[ ] $F_{30} = X^{30} + 522 (X^{25} Y^5 - X^5 Y^{25}) - 10005 (X^{20} Y^{10}+ X^{10} Y^{20}) + Y^{30}$, 
\item[ ] $F_{20} = X^{20} - 228 (X^{15} Y^5 - X^5 Y^{15}) + 494 X^{10} Y^{10} + Y^{20}$ \, and
\item[ ] $F_{12} = XY (X^{10} + 11 X^5 Y^5 -Y^{10})$. 
\end{namelist}

For these polynomials, we have well-known facts as follows.

\begin{fact}
Let $\zeta_5$ be a primitive $5$th root of unity and 
put  
$$\alpha = - \left(\begin{array}{cc}
\zeta_5^3 & 0          \\
0         & \zeta_5^2  \\
\end{array}\right), \quad 
\beta  = \left(\begin{array}{cc}
 0 & 1  \\
-1 & 0  \\
\end{array}\right), \quad $$
$$\gamma   = \frac{1}{\zeta_5^2 - \zeta_5^3} \left(\begin{array}{cc}
 \zeta_5 + \zeta_5^{-1} & 1                          \\
 1                      & -(\zeta_5 + \zeta_5^{-1})  \\
\end{array}\right)$$
and $\overline{I} = \langle \alpha, \beta, \gamma \rangle$. 
Then 
$\mathbb{C} [X,Y]^{\overline{I}} = \mathbb{C}[F_{30}, F_{20}, F_{12}]$. 
Note that $\overline{I} \cong {\rm SL}(2, 5)$: the binary icosahedral subgroup of ${\rm SL}(2, {\mathbb C})$. 
\end{fact}

Under the situation above, 
our main results are stated as follows.

\begin{theorem}\label{main}
Let $C_{30}$, $C_{20}$ and $C_{12}$ be the plane curves defined by

\begin{namelist}{111}
\item[] $C_{30} \colon Z^{30} + F_{30} (X, Y) = 0$, 
\item[] $C_{20} \colon Z^{20} + F_{20} (X, Y) = 0$ \, and  
\item[] $C_{12} \colon Z^{12} + F_{12} (X, Y) = 0$. 
\end{namelist}
Then $|{\rm Aut}(C_d)| = 60 d$ $($$d = 30$, $20$, $12$$)$. 
Furthermore, the following hold:

\begin{namelist}{111}
\item[] ${\rm Aut}(C_{30}) \cong {\mathbb Z}_{15} \times {\rm SL}(2, 5)$, 
\item[] ${\rm Aut}(C_{20}) \cong {\mathbb Z}_5 \times ( {\rm SL}(2, 5) \rtimes {\mathbb Z}_2)$ \, and 
\item[] ${\rm Aut}(C_{12}) \cong {\mathbb Z}_3 \times ( {\rm SL}(2, 5) \rtimes {\mathbb Z}_2)$.
\end{namelist}
\end{theorem}

%%%%%%%%%%%%%%%%%%%%%%%%%%%%%%%%%%%%%%%%%%%%%%%%%%%%%%%%%
\section{Proofs of Theorem \ref{main}}
First of all, we review Theorem \ref{harui_th} {\rm (a-ii)} from the viewpoint of Galois points. 
Let $C$ be a smooth plane curve of degree $d \geq 4$ with a unique Galois point $P$, 
$G$ a subgroup of ${\rm Aut}(C)$. 
Then by Proposition \ref{stform}, 
we may assume that the defining equation of $C$ is given by $Z^d + F_d(X, Y) = 0$ 
for some homogeneous polynomial $F_d (X, Y)$ of degree $d$ and $P = (0:0:1)$. 
Let $\pi_ P : {\mathbb P}^2 \cdots \rightarrow {\mathbb P}^1$ be the projection with the center $P$. 
Then $\pi_P$ is represented as $\pi_P ((X:Y:Z)) = (X:Y)$. 
The Galois group $G_P$ is represented by 
$$G_P = \left\langle \left(\begin{array}{ccc}
1       & 0       & 0  \\
0       & 1       & 0  \\
0       & 0       & \zeta_{d}
\end{array}\right) \right\rangle ,$$
where $\zeta_{d}$ is a primitive $d$-th root of unity.  
We denote by $\lambda_d$ this matrix generating $G_P$. 
Then we get the following commutative diagram as in Theorem \ref{harui_th} (a-ii): 
$$\xymatrix{
1 \ar[r] &  {\mathbb C}^{\times} \ar[r]  & {\rm PBD}(2, 1) \ar[r]^{\rho} & {\rm PGL} (2, {\mathbb C}) \ar[r] & 1 \\
1 \ar[r] &  N       \ar@{^{(}->}[u] \ar[r] & G        \ar@{^{(}->}[u] \ar[r] &  G' \ar@{^{(}->}[u] \ar[r]             & 1. \\
}$$
In this case $N = G_P$. Thus we get the exact sequence  
$$
(\sharp) \hspace{12pt} 1
\rightarrow G_P 
\rightarrow G
\overset{\rho}{\rightarrow} G' 
\rightarrow 1, 
$$
where $G' \subset {\rm PGL} (2, {\mathbb C})$.

Now, we put
$$\sigma = \left(\begin{array}{ccc}
\frac{\zeta_5 - \zeta_5^4}{\sqrt{5}}   & \frac{\zeta_5^3 - \zeta_5^2}{\sqrt{5}} & 0  \\
\frac{\zeta_5^3 - \zeta_5^2}{\sqrt{5}} & - \frac{\zeta_5 - \zeta_5^4}{\sqrt{5}}  & 0  \\
0                                      & 0                                      & 1
\end{array}\right), \quad 
\tau = \left(\begin{array}{ccc}
\zeta_5 & 0 & 0  \\
0       & 1 & 0  \\
0       & 0 & 1
\end{array}\right),$$

$$\rho = \left(\begin{array}{ccc}
\zeta_5^3 & 0       & 0  \\
0       & \zeta_5^2 & 0  \\
0       & 0         & 1
\end{array}\right) \quad {\rm and} \quad 
\varphi  = \left(\begin{array}{ccc}
\zeta_5 & 0 & 0  \\
0       & 1 & 0  \\
0       & 0 & \xi
\end{array}\right),$$
where $\xi^{12} = \zeta_5$. 
We also put  

$$\alpha' = \left(\begin{array}{@{\,}ccc|c@{\,}}
        &          &           &       \\
        & \mbox{\smash{\huge\textit{$\alpha$}}}   &           &       \\
        &          &           &       \\ \hline
        &          &           &    1  \\        
\end{array}\right), 
\beta' = \left(\begin{array}{@{\,}ccc|c@{\,}}
        &          &           &       \\
        & \mbox{\smash{\huge\textit{$\beta$}}}    &           &       \\
        &          &           &       \\ \hline
        &          &           &    1  \\        
\end{array}\right), 
\gamma' = \left(\begin{array}{@{\,}ccc|c@{\,}}
        &          &           &       \\
        & \mbox{\smash{\huge\textit{$\gamma$}}}   &           &       \\
        &          &           &       \\ \hline
        &          &           &    1  \\        
\end{array}\right) \in {\rm GL}(3, {\mathbb C}). $$

Referring to \cite{B}, we see that the image of $\overline{I}$ under the natural homomorphism 
${\rm SL}(2,{\mathbb C}) \rightarrow {\rm PGL}(2,{\mathbb C})$ is isomorphic to $A_5$. 
Further, we define ${\rm S}(2,1) : = \left\langle \alpha', \beta', \gamma' \right\rangle \cong \overline{I}$.

First we deal with $C_{30}$. 
Put $\widetilde{G_0} = \left\langle \sigma, \tau, \lambda_{30} \right\rangle \subset {\rm GL (3, \mathbb{C})}$ 
and $H = \left\langle \sigma, \tau \right\rangle$. 
Then we can check that $(\sigma \tau^4)^2 = 
\left(\begin{array}{ccc}
\zeta_5 & 0       & 0  \\
0       & \zeta_5 & 0  \\
0       & 0       & 1
\end{array}\right)$ and $\rho = \tau \left(\begin{array}{ccc}
\zeta_5 & 0       & 0  \\
0       & \zeta_5 & 0  \\
0       & 0       & 1
\end{array}\right)^2$. 
So we have $\rho \in H$. 
Furthermore, since $\alpha' = (\rho^2 \sigma \rho)^2 \rho$, 
$\beta' = \sigma (\rho^2 \sigma \rho)^2 \rho$ and 
$\gamma' = \sigma^2 \rho$, 
we obtain $H \supset {\rm S}(2,1)$.

We also remark that $\left(\begin{array}{ccc}
-1      & 0       & 0  \\
0       & -1      & 0  \\
0       & 0       & 1
\end{array}\right) \in {\rm S}(2,1)$ and 
$\sigma = \left(\begin{array}{ccc}
-1      & 0       & 0  \\
0       & -1      & 0  \\
0       & 0       & 1
\end{array}\right) \beta' \gamma'$, 
$\tau = \left(\begin{array}{ccc}
-1      & 0       & 0  \\
0       & -1      & 0  \\
0       & 0       & 1
\end{array}\right) \alpha' 
\left(\begin{array}{ccc}
\zeta_5      & 0       & 0  \\
0            & \zeta_5 & 0  \\
0            & 0       & 1
\end{array}\right)^2$. 
Thus we obtain $\left\langle {\rm S}(2,1), 
\left(\begin{array}{ccc}
\zeta_5      & 0       & 0  \\
0            & \zeta_5 & 0  \\
0            & 0       & 1
\end{array}\right)
\right\rangle = 
{\rm S}(2,1) \times
\left\langle
\left(\begin{array}{ccc}
\zeta_5      & 0       & 0  \\
0            & \zeta_5 & 0  \\
0            & 0       & 1
\end{array}\right)
\right\rangle
= H$.

Therefore, we have that
$\widetilde{G_0} = H \times \left\langle
\left(\begin{array}{ccc}
1            & 0       & 0  \\
0            & 1       & 0  \\
0            & 0       & \zeta_{30}
\end{array}\right)
\right\rangle = 
{\rm S}(2,1) \times \left\langle \left(\begin{array}{ccc}
\zeta_5      & 0       & 0  \\
0            & \zeta_5 & 0  \\
0            & 0       & 1
\end{array}\right)
\right\rangle \times \left\langle
\left(\begin{array}{ccc}
1            & 0       & 0  \\
0            & 1       & 0  \\
0            & 0       & \zeta_{30}
\end{array}\right) \right\rangle$.

Put 
$Z := \widetilde{G_0} \cap \left\{ 
\left. \left(\begin{array}{ccc}
\eta         & 0       & 0  \\
0            & \eta    & 0  \\
0            & 0       & \eta
\end{array}\right) \right| \, \eta \in {\mathbb C}^{\times}
\right\} = \left\langle \left(\begin{array}{ccc}
-1           & 0       & 0  \\
0            & -1      & 0  \\
0            & 0       & -1
\end{array}\right), 
\left(\begin{array}{ccc}
\zeta_5      & 0            & 0  \\
0            & \zeta_5      & 0  \\
0            & 0            & \zeta_5
\end{array}\right)
\right\rangle .$ 
Then $G_0 := \widetilde{G_0} / Z \subset G$ and 
$$G_0 = \frac{{\rm S}(2,1) \times \left\langle \left(\begin{array}{ccc}
1            & 0       & 0  \\
0            & 1       & 0  \\
0            & 0       & -1
\end{array}\right)
\right\rangle}
{\left\langle \left(\begin{array}{ccc}
-1           & 0       & 0  \\
0            & -1      & 0  \\
0            & 0       & -1
\end{array}\right)\right\rangle} \times 
\frac{\left\langle \left(\begin{array}{ccc}
\zeta_5      & 0            & 0  \\
0            & \zeta_5      & 0  \\
0            & 0            & 1
\end{array}\right)
\right\rangle \times \left\langle
\left(\begin{array}{ccc}
1            & 0       & 0  \\
0            & 1       & 0  \\
0            & 0       & \zeta_{15}
\end{array}\right)
\right\rangle}{\left\langle\left(\begin{array}{ccc}
\zeta_5      & 0            & 0  \\
0            & \zeta_5      & 0  \\
0            & 0            & \zeta_5
\end{array}\right)
\right\rangle} . $$

Hence $G_0 \cong {\rm SL}(2, 5) \times {\mathbb Z}_{15}$. 
In particular, $\vert G_0 \vert = 120 \cdot 15 = 1800$. 
On the other hand, we see that $\vert G \vert = 30 \cdot 60 = 1800$ by $(\sharp )$. 
Hence $G_0 = G$, which completes the proof of this case.

By a similar argument to the above, we can prove the other cases. 
So, we give the proofs in brief.

For the curve $C_{20}$, we put 
$\widetilde{G_0} = \left\langle \sigma, \tau, \lambda_{20} \right\rangle 
\subset {\rm GL (3, \mathbb{C})}$. 
We see that 
$\widetilde{G_0} = 
{\rm S}(2,1) \times \left\langle \left(\begin{array}{ccc}
1            & 0       & 0  \\
0            & 1       & 0  \\
0            & 0       & \sqrt{-1}
\end{array}\right)
\right\rangle 
\times \left\langle \left(\begin{array}{ccc}
\zeta_5      & 0       & 0  \\
0            & \zeta_5 & 0  \\
0            & 0       & 1
\end{array}\right)
\right\rangle \times \left\langle
\left(\begin{array}{ccc}
1            & 0       & 0  \\
0            & 1       & 0  \\
0            & 0       & \zeta_{5}
\end{array}\right) \right\rangle$. 
Since its center $Z$ is $\left\langle \left(\begin{array}{ccc}
-1           & 0       & 0  \\
0            & -1      & 0  \\
0            & 0       & -1
\end{array}\right), 
\left(\begin{array}{ccc}
\zeta_5      & 0            & 0  \\
0            & \zeta_5      & 0  \\
0            & 0            & \zeta_5
\end{array}\right)
\right\rangle$, 
we obtain 
{\small
$$G_0  = \frac{{\rm S}(2,1) \times \left\langle \left(\begin{array}{ccc}
1            & 0       & 0  \\
0            & 1       & 0  \\
0            & 0       & \sqrt{-1}
\end{array}\right)
\right\rangle}
{\left\langle \left(\begin{array}{ccc}
-1           & 0       & 0  \\
0            & -1      & 0  \\
0            & 0       & -1
\end{array}\right)\right\rangle} \times 
\frac{\left\langle \left(\begin{array}{ccc}
\zeta_5      & 0            & 0  \\
0            & \zeta_5      & 0  \\
0            & 0            & 1
\end{array}\right)
\right\rangle \times \left\langle
\left(\begin{array}{ccc}
1            & 0       & 0  \\
0            & 1       & 0  \\
0            & 0       & \zeta_5
\end{array}\right)
\right\rangle}{\left\langle\left(\begin{array}{ccc}
\zeta_5      & 0            & 0  \\
0            & \zeta_5      & 0  \\
0            & 0            & \zeta_5
\end{array}\right)
\right\rangle}, $$
}
where $G_0 = \widetilde{G_0} / Z \subset G$.

Further, since ${\rm SL}(2,5) \cong \dfrac{{\rm S}(2,1) \times \left\langle \left(\begin{array}{ccc}
1            & 0       & 0  \\
0            & 1       & 0  \\
0            & 0       & -1
\end{array}\right)
\right\rangle}
{\left\langle \left(\begin{array}{ccc}
-1           & 0       & 0  \\
0            & -1      & 0  \\
0            & 0       & -1
\end{array}\right)\right\rangle}$, 
we have the following exact sequence: 
$$
1 \rightarrow 
{\rm SL}(2,5) \rightarrow 
\frac{{\rm S}(2,1) \times \left\langle \left(\begin{array}{ccc}
1            & 0       & 0  \\
0            & 1       & 0  \\
0            & 0       & \sqrt{-1}
\end{array}\right)
\right\rangle}
{\left\langle \left(\begin{array}{ccc}
-1           & 0       & 0  \\
0            & -1      & 0  \\
0            & 0       & -1
\end{array}\right)\right\rangle}
\overset{\delta}{\rightarrow}
\left\{  \pm 1 \right\} \rightarrow
1,
$$
where $\delta : \left(\begin{array}{@{\,}ccc|c@{\,}}
        &          &                       &     \\
        & \mbox{\smash{\huge\textit{A}}}   &  &  \\
        &          &                       &     \\ \hline
        &          &                       &     \alpha  \\        
\end{array}\right) \mapsto \alpha^2$. 
The sequence is split by 
$-1 \mapsto \left(\begin{array}{ccc}
0            & 1       & 0  \\
-1           & 0       & 0  \\
0            & 0       & \sqrt{-1}
\end{array}\right)$.

Hence $G_0 \cong  ( {\rm SL}(2, 5) \rtimes {\mathbb Z}_2) \times {\mathbb Z}_5$. 
In particular, $\vert G_0 \vert = 120 \cdot 2 \cdot 5 = 1200$. 
On the other hand, we see that $\vert G \vert = 20 \cdot 60 = 1200$ by $(\sharp )$. 
Hence $G_0 = G$, which completes the proof of this case.

Finally, for the curve $C_{12}$, 
we put 
$\widetilde{G_0} = \left\langle \sigma, \varphi \right\rangle$, and 
$K = \left\langle \sigma, \rho  \right\rangle$. 
We can check that $K = {\rm S}(2, 1)$. 
Putting $\varepsilon := (\sigma \varphi^4)^3 = \left(\begin{array}{ccc}
\zeta_5      & 0            & 0  \\
0            & \zeta_5      & 0  \\
0            & 0            & \zeta_5
\end{array}\right)$, 
we obtain $\widetilde{G_0} = \left\langle K, \lambda_{12}, \varepsilon \right\rangle$. 
Furthermore we get $\widetilde{G_0} = \left\langle K \right\rangle \times \left\langle \lambda_{12} \right\rangle \times  \left\langle \varepsilon \right\rangle$
$$= {\rm S}(2, 1) \times \left\langle \left(\begin{array}{ccc}
 1           & 0       & 0  \\
0            &  1      & 0  \\
0            & 0       & \sqrt{-1}
\end{array}\right)\right\rangle
\times 
\left\langle \left(\begin{array}{ccc}
1            & 0       & 0  \\
0            & 1       & 0  \\
0            & 0       & \omega 
\end{array}\right)\right\rangle
\times 
\left\langle \left(\begin{array}{ccc}
\zeta_5      & 0            & 0  \\
0            & \zeta_5      & 0  \\
0            & 0            & \zeta_5
\end{array}\right)\right\rangle ,$$
where $\omega$ is a cubic root of unity.

Hence $G_0 \cong   ( {\rm SL}(2, 5) \rtimes {\mathbb Z}_2) \times {\mathbb Z}_3$. 
In particular, $\vert G_0 \vert = 120 \cdot 2 \cdot 3 = 720$. 
On the other hand, we see that $\vert G \vert = 12 \cdot 60 = 720$ by $(\sharp )$. 
Hence $G_0 = G$, which completes the proof of this case.

%%%%%%%%%%%%%%%%%%%%%%%%%%%%%%%%%%%%%%%%%%%%%%%%%%%%%%%%%

\end{document}